\magnification=1200

\hsize=11.25cm    
\vsize=18cm     
\parindent=12pt   \parskip=5pt     

\hoffset=.5cm   
\voffset=.8cm   

\pretolerance=500 \tolerance=1000  \brokenpenalty=5000

\catcode`\@=11

\font\eightrm=cmr8         \font\eighti=cmmi8
\font\eightsy=cmsy8        \font\eightbf=cmbx8
\font\eighttt=cmtt8        \font\eightit=cmti8
\font\eightsl=cmsl8        \font\sixrm=cmr6
\font\sixi=cmmi6           \font\sixsy=cmsy6
\font\sixbf=cmbx6

\font\tengoth=eufm10 
\font\eightgoth=eufm8  
\font\sevengoth=eufm7      
\font\sixgoth=eufm6        \font\fivegoth=eufm5

\skewchar\eighti='177 \skewchar\sixi='177
\skewchar\eightsy='60 \skewchar\sixsy='60

\newfam\gothfam           \newfam\bboardfam

\def\tenpoint{
  \textfont0=\tenrm \scriptfont0=\sevenrm \scriptscriptfont0=\fiverm
  \def\rm{\fam\z@\tenrm}
  \textfont1=\teni  \scriptfont1=\seveni  \scriptscriptfont1=\fivei
  \def\oldstyle{\fam\@ne\teni}\let\old=\oldstyle
  \textfont2=\tensy \scriptfont2=\sevensy \scriptscriptfont2=\fivesy
  \textfont\gothfam=\tengoth \scriptfont\gothfam=\sevengoth
  \scriptscriptfont\gothfam=\fivegoth
  \def\goth{\fam\gothfam\tengoth}
  
  \textfont\itfam=\tenit
  \def\it{\fam\itfam\tenit}
  \textfont\slfam=\tensl
  \def\sl{\fam\slfam\tensl}
  \textfont\bffam=\tenbf \scriptfont\bffam=\sevenbf
  \scriptscriptfont\bffam=\fivebf
  \def\bf{\fam\bffam\tenbf}
  \textfont\ttfam=\tentt
  \def\tt{\fam\ttfam\tentt}
  \abovedisplayskip=12pt plus 3pt minus 9pt
  \belowdisplayskip=\abovedisplayskip
  \abovedisplayshortskip=0pt plus 3pt
  \belowdisplayshortskip=4pt plus 3pt 
  \smallskipamount=3pt plus 1pt minus 1pt
  \medskipamount=6pt plus 2pt minus 2pt
  \bigskipamount=12pt plus 4pt minus 4pt
  \normalbaselineskip=12pt
  \setbox\strutbox=\hbox{\vrule height8.5pt depth3.5pt width0pt}
  \let\bigf@nt=\tenrm       \let\smallf@nt=\sevenrm
  \normalbaselines\rm}

\def\eightpoint{
  \textfont0=\eightrm \scriptfont0=\sixrm \scriptscriptfont0=\fiverm
  \def\rm{\fam\z@\eightrm}
  \textfont1=\eighti  \scriptfont1=\sixi  \scriptscriptfont1=\fivei
  \def\oldstyle{\fam\@ne\eighti}\let\old=\oldstyle
  \textfont2=\eightsy \scriptfont2=\sixsy \scriptscriptfont2=\fivesy
  \textfont\gothfam=\eightgoth \scriptfont\gothfam=\sixgoth
  \scriptscriptfont\gothfam=\fivegoth
  \def\goth{\fam\gothfam\eightgoth}
  
  \textfont\itfam=\eightit
  \def\it{\fam\itfam\eightit}
  \textfont\slfam=\eightsl
  \def\sl{\fam\slfam\eightsl}
  \textfont\bffam=\eightbf \scriptfont\bffam=\sixbf
  \scriptscriptfont\bffam=\fivebf
  \def\bf{\fam\bffam\eightbf}
  \textfont\ttfam=\eighttt
  \def\tt{\fam\ttfam\eighttt}
  \abovedisplayskip=9pt plus 3pt minus 9pt
  \belowdisplayskip=\abovedisplayskip
  \abovedisplayshortskip=0pt plus 3pt
  \belowdisplayshortskip=3pt plus 3pt 
  \smallskipamount=2pt plus 1pt minus 1pt
  \medskipamount=4pt plus 2pt minus 1pt
  \bigskipamount=9pt plus 3pt minus 3pt
  \normalbaselineskip=9pt
  \setbox\strutbox=\hbox{\vrule height7pt depth2pt width0pt}
  \let\bigf@nt=\eightrm     \let\smallf@nt=\sixrm
  \normalbaselines\rm}

\tenpoint

\def\pc#1{\bigf@nt#1\smallf@nt}         \def\pd#1 {{\pc#1} }

\frenchspacing

\def\raggedbottom{\topskip 10pt plus 36pt\r@ggedbottomtrue}

\def\pointir{\unskip . --- \ignorespaces}

\def\Medbreak{\vskip-\lastskip\medbreak}

\long\def\th#1 #2\enonce#3\endth{
   \Medbreak\noindent
   {\pc#1} {#2\unskip}\pointir{\it #3}\smallskip}

\def\decale#1{\smallbreak\hskip 28pt\llap{#1}\kern 5pt}
\def\decaledecale#1{\smallbreak\hskip 34pt\llap{#1}\kern 5pt}
\def\puce{\smallbreak\hskip 6pt{$\scriptstyle\bullet$}\kern 5pt}

\def\eqalign#1{\null\,\vcenter{\openup\jot\m@th\ialign{
\strut\hfil$\displaystyle{##}$&$\displaystyle{{}##}$\hfil
&&\quad\strut\hfil$\displaystyle{##}$&$\displaystyle{{}##}$\hfil
\crcr#1\crcr}}\,}

\catcode`\@=12

\showboxbreadth=-1  \showboxdepth=-1

\newcount\numerodesection \numerodesection=1
\def\section#1{\bigbreak
 {\bf\number\numerodesection.\ \ #1}\nobreak\medskip
 \advance\numerodesection by1}

\mathcode`A="7041 \mathcode`B="7042 \mathcode`C="7043 \mathcode`D="7044
\mathcode`E="7045 \mathcode`F="7046 \mathcode`G="7047 \mathcode`H="7048
\mathcode`I="7049 \mathcode`J="704A \mathcode`K="704B \mathcode`L="704C
\mathcode`M="704D \mathcode`N="704E \mathcode`O="704F \mathcode`P="7050
\mathcode`Q="7051 \mathcode`R="7052 \mathcode`S="7053 \mathcode`T="7054
\mathcode`U="7055 \mathcode`V="7056 \mathcode`W="7057 \mathcode`X="7058
\mathcode`Y="7059 \mathcode`Z="705A


\def\hfl#1#2#3{\smash{\mathop{\hbox to#3{\rightarrowfill}}\limits
^{\textstyle#1}_{\textstyle#2}}}

\def\P{{\bf P}}
\def\Q{{\bf Q}}

\def\Z{{\bf Z}}

\def\F{{\bf F}}

\def\GL{{\bf GL}}
\def\AGL{{\bf AGL}}

\def\Hom{\mathop{\rm Hom}\nolimits}

\def\Ext{\mathop{\rm Ext}\nolimits}

\def\Card{\mathop{\rm Card}\nolimits}
\def\Gal{\mathop{\rm Gal}\nolimits}
\def\Ker{\mathop{\rm Ker}\nolimits}

\def\to{\rightarrow}

\def\boxit#1{\vbox{\hrule\hbox{\vrule\kern1pt
       \vbox{\kern1pt#1\kern1pt}\kern1pt\vrule}\hrule}}
\def\cqfd{\hfill\boxit{\phantom{\i}}}

\newcount\numero \numero=1
\def\numeroter{{({\oldstyle\number\numero})}\ \advance\numero by1}

\newcount\refno 
\long\def\ref#1:#2<#3>{                                        
\global\advance\refno by1\par\noindent                              
\llap{[{\bf\number\refno}]\ }{#1} \pointir{\it #2} #3\goodbreak }

\newcount\refno 
\long\def\ref#1:#2<#3>{                                        
\global\advance\refno by1\par\noindent                              
\llap{[{\bf\number\refno}]\ }{#1} \pointir{\it #2} #3\goodbreak }

\def\citer#1(#2){[{\bf\number#1}\if#2\empty\relax\else,\ {#2}\fi]}

\def\boxit#1{\vbox{\hrule\hbox{\vrule\kern1pt
       \vbox{\kern1pt#1\kern1pt}\kern1pt\vrule}\hrule}}
\def\cqfd{\hfill\boxit{\phantom{\i}}}

\openup.25\jot

\newbox\bibbox
\setbox\bibbox\vbox{\bigbreak
\centerline{{\pc BIBLIOGRAPHY}}
\nobreak

\ref{\pc CLIFFORD} (A):
Representations induced in an invariant subgroup,
<Annals of Math.\  {\bf 38} (1937) 3, 533--550.>
\newcount\clifford \global\clifford=\refno

\ref{\pc COX} (D): Galois theory, <Second edition. John Wiley \&
Sons, Inc., Hoboken, NJ, 2012. xxviii+570 pp.> \newcount\cox
\global\cox=\refno

\ref{\pc DALAWAT} (C):
Serre's ``\thinspace formule de masse\thinspace'' in prime degree,
<Monats\-hefte Math.\ {\bf 166} (2012) 1, 73--92.
Cf.~arXiv\string:1004.2016v6.>
\newcount\monatshefte \global\monatshefte=\refno

\ref{\pc DALAWAT} (C) \& {\pc LEE} (JJ):
Tame ramification and group cohomology,
<J.\ Ramanujan Math.\ Soc.\ {\bf 32} (2017) 1, 51--74.
Cf.~arXiv\string:1305.2580v4.>
\newcount\dalawatlee \global\dalawatlee=\refno

\ref{\pc DALAWAT} (C):
$\F_p$-representations over $p$-fields,
<arXiv\string:1608.04181.>
\newcount\irrepp \global\irrepp=\refno

\ref{\pc DALAWAT} (C):
Little galoisian modules,
<arXiv\string:1608.04182.>
\newcount\littlegalmods \global\littlegalmods=\refno

\ref{\pc DALAWAT} (C):
Wildly primitive extensions,
<arXiv\string:1608.04183.>
\newcount\wildprim \global\wildprim=\refno

\ref{\pc DEL \pc CORSO} (I), {\pc DVORNICICH} (R) \& {\pc MONGE} (M):
On wild extensions of a $p$-adic field,
<J.\ Number Theory {\bf 174} (2017), 322--342.  Cf.~aXiv\string:1601.05939.> 
\newcount\deldvomonge \global\deldvomonge=\refno

\ref{\pc FONCTION} (L):
The\/ $L$-functions and modular forms database, 
<online database, http:/$\!$/www.lmfdb.org/>
\newcount\lmfdb \global\lmfdb=\refno

\ref{\pc HENNIART} (G):
Representations du groupe de Weil d'un corps local,
<http:/$\!$/sites.mathdoc.fr/PMO/feuilleter.php?id=PMO\_1979>
\newcount\henniart \global\henniart=\refno

\ref{\pc NEUMANN} (P): The concept of primitivity in group theory
and the Second Memoir of Galois, <Arch.\ Hist.\ Exact Sci.\ {\bf 60}
(2006) 4, 379--429.>
\newcount\neumann \global\neumann=\refno

\ref{\pc OVERFLOW} (M):
Solvable irreducible subgroups of the\/ $\GL_n$ of\/ $\F_p$ ($p$ prime),
<http:/$\!$/mathoverflow.net/q/241982>
\newcount\csdmo \global\csdmo=\refno

} 

\centerline{\bf Solvable primitive extensions} 
\bigskip\bigskip 
\centerline{Chandan Singh Dalawat} 
\centerline{Harish-Chandra Research Institute}
\centerline{Chhatnag Road, Jhunsi, Allahabad 211019, India} 
\centerline{\tt dalawat@gmail.com}

\bigskip\bigskip

{{\bf Abstract}.  A finite separable extension $E$ of a field $F$ is
  called primitive if there are no intermediate extensions.  It is
  called a solvable extension if the group $\Gal(\hat E|F)$ of
  $F$-automorphisms of its galoisian closure $\hat E$ over $F$ is
  solvable. We show that a solvable primitive extension $E$ of
  $F$ is uniquely determined (up to $F$-isomorphism) by $\hat E$ and
  characterise the extensions $D$ of $F$ such that $D=\hat E$ for some
  solvable primitive extension $E$ of $F$.

\footnote{}{{\it MSC2010~:} Primary 12F10, 12G05, Seconday 20B15, 20C05}
\footnote{}{{\it Keywords~:} Solvable extensions, primitive extensions, $p$-extensions
}}

\bigskip\bigskip

\numeroter Let $F$ be a field.  All extensions of $F$ appearing below
are assumed to be {\it separable\/} over $F$.  A finite extension $E$
of $F$ of degree $[E:F]>1$ is called {\it primitive\/} if the only
intermediate fields $F\subset K\subset E$ are $K=F$ and $K=E$.  We say
that $E$ is {\it solvable\/} over $F$ if the group $G=\Gal(\hat E|F)$
is solvable, where $\hat E$ is the galoisian closure of $E$ over $F$.
Galois proved (see below) that if $E$ is a solvable primitive
extension of $F$, then $[E:F]=l^n$ for some prime $l$ and some $n>0$.
So we will later fix the prime $l$ and the integer $n$ as well~;
extensions whose degree is a power of~$l$ will simply be called
$l$-{\it extensions}.  The problem is to {\it parametrise all solvable
primitive extensions of degree\/~$l^n$ over\/ $F$}.

\numeroter  Let $\Omega$ be a finite set with $\Card(\Omega)>1$.  A
partition $(\Omega_i)_{i\in I}$ of $\Omega$ will be called {\it
essential\/} if $\Card(I)>1$, $\Card(\Omega_i)>0$ for all $i\in I$ and
$\Card(\Omega_i)>1$ for some $i\in I$.  Let $G$ be a subgroup of the
symmetric group ${\goth S}_\Omega$. A partition $(\Omega_i)_{i\in I}$
of $\Omega$ will be called $G$-{\it stable\/} if $G$ permutes the
parts~: if there is a map $\pi:G\to{\goth S}_I$ such that
$g.\Omega_i=\Omega_{\pi(g)(i)}$ for every $g\in G$ and every $i\in
I$~; if so, $\pi$ is a homomorphism of groups.  A subgroup
$G\subset{\goth S}_\Omega$ is called {\it primitive\/} if its order is
$>1$, and no $G$-stable essential partition of $\Omega$ exists.  Every
primitive subgroup of ${\goth S}_\Omega$ is transitive.

\numeroter (Galois) {\it If\/ $G$ is a solvable primitive subgroup
of\/ ${\goth S}_\Omega$, then there is a unique structure on\/
$\Omega$ of an affine $n$-space over $\F_l$ (for some prime\/ $l$ and
some\/ $n>0$) such that
$$
N\subset G\subset\AGL(\Omega)\subset{\goth S}_\Omega,
$$
where\/ $N$ (resp.~$\AGL(\Omega)$) is the space of translations\/
(resp.~the group of automorphisms or invertible affine maps) of the
affine space\/ $\Omega$}.

{\it Proof}. The group $G$ is not trivial by hypothesis.  Let $N$ be a
minimal normal subgroup of $G$.  Since $G$ is solvable, $N$ is a
vector $n$-space over $\F_l$ for some prime\/ $l$ and some\/ $n>0$.
Since $G$ is primitive, $N$ is transitive.  Since $N$ is commutative
and transitive, $\Omega$ is an $N$-torsor (an affine space over $\F_l$
whose space of translations is $N$).  Finally, one checks that $G$
acts on $\Omega$ by affine maps (because the conjugation action of
$G/N$ on $N$ is by $\F_l$-linear maps, automatically).  For more
details and some historical remarks, see \citer\cox(p.~435)
or \citer\neumann(p.~420). \cqfd

\numeroter {\it Let\/ $\Omega$ be an affine space over\/ $\F_l$ of
dimension\/ $n>0$ and let\/ $N$ be its space of translations.  An
intermediate group\/ $N\subset G\subset\AGL(\Omega)$ is solvable
(resp.~primitive) if and only if\/ $G/N$ is solvable (resp.~the\/
$\F_l[G/N]$-module\/ $N$ is simple)}.

{\it Proof}. The bit about solvability is clear.  Suppose that $G$ is
imprimitive, and let $(\Omega_i)_{i\in I}$ be a $G$-stable essential
partition of $\Omega$ $\oldstyle(2)$.  Since $G$ is transitive (even
$N$ is transitive), the parts $\Omega_i$ have the same cardinal $l^m$,
for some $m\in[1,n[$.  One checks that they are affine subspaces, all
parallel to each other.  Their common direction $M\subset N$ is a
$G$-stable subspace of dimension~$m$, so $N$ is not simple as an
$\F_l[G/N]$-module.  Conversely, if the $\F_l[G/N]$-module $N$ is not
simple, let $M\subset N$ be a $G$-stable subspace of some dimension
$m\in[1,n[$.  The family of affine subspaces of $\Omega$ of direction
$M$ is a $G$-stable essential partition of $\Omega$, so $G$ is
imprimitive. For more details, see \citer\cox(p.~436).  \cqfd


\numeroter {\it A finite extension\/ $E$ of\/ $F$ is primitive\/
$\oldstyle(1)$ if and only if the group\/ $G=\Gal(\hat E|F)$ is a
primitive subgroup\/ $\oldstyle(2)$ of\/ ${\goth S}_\Omega$, where
$\Omega=\Hom_F(E,\hat E)$}.

{\it Proof}. Note first that $G$ is a transitive subgroup of ${\goth
S}_\Omega$.  Clearly the extension $E$ is primitive over $F$ if and
only if the subgroup $H=\Gal(\hat E|E)$ is maximal in $G$.  Now, $H$
is the stabiliser of $j_E\in\Omega$, where $j_E$ is the inclusion of
$E$ in $\hat E$, and stabilisers of points are maximal subgroups if
and only if the transitive permutation group is primitive.  \cqfd


\goodbreak

\numeroter  {\it If\/ $E$ is a solvable primitive extension of\/ $F$,
then\/ $[E:F]=l^n$ for some prime\/ $l$ and some\/ $n>0$.}

{\it Proof}.  Since $G=\Gal(\hat E|F)$ is a solvable by definition
$\oldstyle(1)$ and primitive $\oldstyle(5)$ as a subgroup of ${\goth
S}_\Omega$ ($\Omega=\Hom_F(E,\hat E)$), one has $\Card\Omega=l^n$
for some prime\/ $l$ and some\/ $n>0$ $\oldstyle(3)$.  But
$[E:F]=\Card\Omega$, hence the result.  \cqfd

\goodbreak 

\numeroter This being so, let us fix the prime $l$ and the integer
$n>0$, and parametrise the set of solvable primitive extensions $E$ of
$F$ of degree~$l^n$.  We will show that $E$ is determined (up to
$F$-isomorphism) by its galoisian closure $\hat E$, and characterise
the galoisian extensions of $F$ which arise in this way, thereby
establishing a canonical bijection between the set of such $E$ and the
set of certain pairs $(\rho,D)$ which we now describe.  Basically, the
pair we attach to $E$ is a way of encapsulating the information
carried by $\hat E$, so let us begin by making that information
explicit.

\numeroter  We have seen ($\oldstyle3$) that the minimal normal
subgroup $N$ of $\Gal(\hat E|F)$ is an $\F_l$-space of dimension~$n$,
and it is a faithful simple module ($\oldstyle4$) over $\Gal(K|F)$,
where $K=\hat E^N$.  Denote by $\rho$ this representation of
$\Gal(K|F)$ on $N$.  As $\Gal(K|F)$ is a quotient of $\Gal(\tilde
F|F)$, where $\tilde F$ is a maximal galoisian extension of the field
$F$, $\rho$ can be viewed as a representation of $\Gal(\tilde F|F)$~;
we then have $K=\tilde F^{\Ker\rho}$.  So from $E$ we get the pair
$(\rho,\hat E)$ consisting of an irreducible $\F_l$-representation
$\rho$ of $\Gal(\tilde F|F)$ of degree~$n$ with solvable image, and an
elementary abelian $l$-extension (namely $\hat E$) of the fixed field
$K$ of the kernel of $\rho$ which is galoisian over $F$ and such that
the conjugation action of $\Gal(K|F)$ on $N=\Gal(\hat E|K)$ resulting
from the short exact sequence $1\to N\to\Gal(\hat E|
F)\to\Gal(K|F)\to1$ is given by $\rho$.

\numeroter  Conversely, suppose that we are given a pair $(\rho,D)$
consisting of an irreducible representation $\rho$ of $\Gal(\tilde
F|F)$ on some $\F_l$-space $N$ of dimension~$n$ with solvable image,
and an $N$-extension $D$ of the fixed field $K$ of the kernel of
$\rho$ such that $D$ is galoisian over $F$ and such that the resulting
conjugation action of $\Gal(K|F)$ on $N=\Gal(D|K)$ is given by $\rho$.
We will show that such a pair $(\rho,D)$ determines a solvable
primitive extension $E$ of $F$ of degree $l^n$, unique up to
$F$-isomorphism, such that the pair associated to $E$ by the
construction ($\oldstyle8$) is $(\rho,D)$.  This will establish the
desired canonical bijection between the set of such $E$ and the set of
such pairs $(\rho,D)$.  The point of having such a bijection is that
for certain fields $F$, the pairs $(\rho,D)$ can be explicitly
computed, as is explained briefly in ($\oldstyle18$).  Let's get
going.

\numeroter  {\it Let\/ $G$ be a finite group and\/ $M$ an\/
$\F_l[G]$-module.  If\/ $G$ is solvable and\/ $M$ is faithful and
simple, then\/ $H^i(G,M)=\{0\}$ for every\/ $i>0$.}

{\it Proof\/} (Jeremy Rickard \citer\csdmo()).  This is clear if the
group $G$ is trivial (in which case $M=\F_l$), so assume that the
order of $G$ is $>1$ (so that $M\neq\F_l$).  Let $N$ be a minimal
normal sugroup of $G$.  Since $G$ is solvable, $N$ is an $\F_p$-space
(for some prime $p$) of some dimension $a>0$.  By Clifford's theorem
(the restriction of an irreducible representation to a normal subgroup
is semisimple) \citer\clifford(), the $\F_l[N]$-module $M$ is
semisimple.

If we had $p=l$, the $\F_l[N]$-module $M$ would be trivial,
contradicting the faithfulness of the $G$-module $M$.  So $p\neq l$.
Consider the central idempotent $r_N=p^{-a}\sum_{n\in N}n$ of
$\F_l[G]$, and put $s_N=1-r_N$.  Clearly, we have a {\it direct sum\/}
decomposition as an $\F_l[G]$-module~: $M=r_NM+s_NM$.  Since $M$ is
simple, we must have $M=r_NM$ or $M=s_NM$.

If we had $M=r_NM$, the action of $N$ on $M$ would be trivial (since
the only simple $\F_l[N]$-module which is not killed by $r_N$ is
$\F_l$), again contradicting the faithfulness of the $G$-module $M$.
Therefore $M=s_NM$, which implies that $M$ and the trivial module
$\F_l$ are in different blocks of $\F_l[G]$ (whatever that may mean),
and hence $H^i(G,M)=\Ext^i_{\F_l[G]}(\F_l,M)=\{0\}$ for every $i>0$.
C'est gagn\'e~!  \cqfd

\smallbreak

\numeroter  {\it Let\/ $G$ be a finite solvable group and let\/ $M$
be a faithful simple\/ $\F_l[G]$-module.  Every extension\/ $L$ of\/
$G$ by the $G$-module\/ $M$ (any short exact sequence\/ $\{1\}\to M\to
L\to G\to\{1\}$ in which the conjugation action of\/ $G$ on\/ $M$ is
the given module structure) splits, and any two sections $G\to L$ are
conjugate in $L$.}

{\it Proof}.  This follows from the case $i=2$ (resp.~$i=1$) of
$\oldstyle(10)$.  Indeed, extensions of $G$ by the $G$-module $M$ are
classified by the group $H^2(G,M)$, and we have seen that this group
is trivial.  Similarly, conjugacy classes of sections of the split
extension of $G$ by $M$ are classified by the group $H^1(G,M)$, and
this group too has been shown to be trivial.  \cqfd

\numeroter {\it Remark}.  See \citer\csdmo() for a direct proof that
in a solvable primitive permutation group $L$ with minimal normal
subgroup $M$ (an $\F_l$-space), there is a unique conjugacy class of
complements to $M$ (all of them maximal subgroups, since each is the
stabiliser of a point).  I thank Peter Neumann for pointing out that
this is covered in the books by Doerk \& Hawkes (p.~55), Huppert I
(p.~159), and Suzuki II (p.~102).


\numeroter {\it Every pair\/ $(\rho,D)$ as in\/ $\oldstyle(9)$
comes, via the construction\/ $\oldstyle(8)$, from a unique solvable
primitive extension\/ $E$ of\/ $F$ of degree\/ $l^n$.}

{\it Proof}.  Let $K=\tilde F^{\Ker\rho}$ be the fixed field of the
kernel of $\rho$.  We have seen that the short exact sequence
$$
1\to\Gal(D|K)\to\Gal(D|F)\to\Gal(K|F)\to1
$$
splits and any two sections are conjugate $\oldstyle(11)$, so there is
an extension $E$ of $F$ (unique up to $F$-isomorphism) linearly
disjoint from $K$ and such that $\hat E=EK=D$.  Since the group
$\Gal(D|F)$ is solvable and primitive $\oldstyle(4)$, and the minimal
normal subgroup $\Gal(D|K)$ is an $\F_l$-space of dimension~$n$, the
extension $E$ of $F$ is solvable and primitive $\oldstyle(5)$ of
degree $l^n$.  \cqfd

\smallbreak

\numeroter {\bf Summary}.  A solvable primitive $l$-extension $E$ of
$F$ is uniquely determined (up to $F$-isomorphism) by its galoisian
closure $\hat E$ over $F$.  A galoisian extension $D$ of $F$ is of the
form $\hat E$ for some solvable primitive $l$-extension $E$ of $F$ if
and only if $D$ is an elementary abelian $l$-extension of the fixed
field $K$ of the kernel of an irreducible $\F_l$-representation $\rho$
of $\Gal(\tilde F|F)$ with solvable image such that $D$ is galoisian
over $F$ and the resulting conjugation action of $\Gal(K|F)$ on
$\Gal(D|K)$ is given by $\rho$.  In terms of the parameter $(\rho,D)$
of $E$, the group $\Gal(\hat E|F)$ is given by
$\Gal(D|K)\times_\rho\Gal(K|F)$.

This parametrisation is useful when $F$ is a $p$-field (a local field
with finite resiude field of characteristic~$p$) because we can
classify irreducible $\F_l$-representations $\rho$ of $\Gal(\tilde
F|F)$ and, for every such $\rho$, determine the structure of the
$\F_l[\Gal(K|F)]$-module $\Gal(M|K)$, where $K$ is the fixed field of
the kernel of $\rho$ and $M$ is the maximal abelian extension of $K$
of exponent~$l$.  See $\oldstyle(18)$ for a brief discussion.

\smallbreak

\numeroter  {\bf Primitive quartic extensions}.  Taking  $l=2$ and
$n=2$, one can parametrise primitive quartic extensions (solvability
is automatic because ${\goth S}_4$ is solvable).  The action of
$\GL_2(\F_2)$ on the set $\P_1(\F_2)$ consisting of the three lines in
$\F_2^2$ identifies it with ${\goth S}_3$.  With this identification,
the only subgroups for which $\F_2^2$ is a simple module are ${\goth
A}_3$ and ${\goth S}_3$~; denote these irreducible representations by
$a_3$ and $s_3$ respectively.  So there are two kinds of primitive
quartic extensions $E$ depending upon whether the group $\Gal(\hat
E|F)$ is isomorphic to ${\goth A}_4=\F_2^2\times_{a_3}{\goth A}_3$ or
to ${\goth S}_4=\F_2^2\times_{s_3}{\goth S}_3$, where $\hat E$ is the
galoisian closure of $E$ over $F$.  Conversely, given any ${\goth
A}_4$- or ${\goth S}_4$-extension $D$, there is a unique (up to
isomorphism) primitive quartic extension $E$ such that $\hat E=D$.

\smallbreak

\numeroter  {\bf Primitive quartic fields}.  When $F=\Q$,
class field theory can be used to determine all such~$D$ and hence all
primitive quartic fields~$E$.  According to the
database \citer\lmfdb(), the smallest (primitive) quartic field $E$
with $\Gal(\hat E|\Q)$ isomorphic to ${\goth A}_4$ (resp.~${\goth
S}_4$) is the one defined by $x^4 - 2x^3 + 2x^2 + 2$ (resp.~$x^4 - x +
1$).  In the first case, the corresponding ${\goth A}_3$-field is the
one defined by $x^3-x^2-2x+1$, and in the second case the ${\goth
S}_3$-field is the galoisian closure of the cubic field defined by
$x^3-4x-1$.  I thank John Jones for confirming this information.

\smallbreak

\numeroter {\bf Primitive quartic extensions of $\Q_2$}.    Let $K$ be
the unique $(\Z/3\Z)$-extension of $\Q_2$.  It can be checked that $K$
  has a unique biquadratic extension $D$ which is galoisian over
  $\Q_2$ and such that the conjugation action of $\Gal(K|\Q_2)=\Z/3\Z$
  on $\Gal(D|K)$ is through $a_3$ $\oldstyle(15)$.  This $D$ is thus
  the unique ${\goth A}_4$-extension of $\Q_2$~; up to isomorphism,
  the corresponding primitive quartic extension is $\Q_2(x)$, where
  $x^4-2x^2+2x-2=0$.

Now let $K$ stand for the unique ${\goth S}_3$-extension of $\Q_2$.
It can be checked that $K$ has {\it three\/} biquadratic extensions
$D$ which are galoisian over $\Q_2$ and such that the conjugation
action of $\Gal(K|\Q_2)={\goth S}_3$ on $\Gal(D|K)$ is through $s_3$
$\oldstyle(15)$.  Thus there are precisely three ${\goth
S}_4$-extensions $D$ of $\Q_2$.  The corresponding primitive quartic
extensions of $\Q_2$ are given (up to isomorphism) by
$$
x^4-2x+2,\ \
x^4-4x+2,\ \
x^4-4x^2+4x-2.
$$
The four polynomials defining the four primitive quartic extensions of
$\Q_2$ are taken from \citer\henniart(p.~111).  Of course, one can
arrive at them directly by examining all polynomials
$x^4+2ax^3+2bx^2+2cx+2d$ with $a,b,c\in\Z_2$ and $d\in\Z_2^\times$
(Eisenstein polynomials).

\numeroter  {\bf The local theory}. The main aim of this series of four
Notes, of which this is the first and purely algebraic one, is to
 parametrise the set of (solvable) primitive $l$-extensions of a
 $p$-field $F$ --- a local field with finite residue field of
 characteristic~$p$.  This problem was already considered by Krasner
 in the 1930s, at least when $F$ is a finite extension of $\Q_p$.
 Every finite (separable) extension of a $p$-field is solvable, and
 the case $l=p$ is the only interesting one.  Thus, our final result
 will be a parametrisation of the set of all primitive $p$-extensions
 of $p$-fields.  Primitive quartic extensions of $\Q_2$ discussed in
 $\oldstyle(17)$ will turn out to be the simplest special case of the
 general theory.

In order to carry out the approach summarised in $\oldstyle(14)$, the
second paper \citer\irrepp() of this series determines the irreducible
$\F_p$-representations $\rho$ of $\Gal(\tilde F|F)$, mostly following
Koch.  We will observe that if we fix~$n>0$, then there is a certain
explicit finite tamely ramified split galoisian extension $L_n$ of $F$
which contains the fixed field of the kernel of every $\rho$ of
degree~$n$ (and such that, if $F$ has characteristic~$0$, then
$L_n^\times$ has an element of order~$p$).

Here, a tamely ramified galoisian extension $L$ of $F$ is said to be
{\it split\/} if, $L_0$ being the maximal unramified extension of $F$
in $L$, the short exact sequence
$$
1\to\Gal(L|L_0)\to\Gal(L|F)\to\Gal(L_0|F)\to1
$$
splits.  It can be shown that the splitting of this sequence is
equivalent to the existence of an intermediate extension $F\subset
E\subset L$ which is totally ramified (but not necessarily galoisian)
over $F$ and such that $[E:F]=[L:L_0]$.
See \citer\dalawatlee(Remark~7.1.4), for example.

The elementary abelian $p$-extensions $D$ of the fixed field $K$ of
the kernel of $\rho$ which are galoisian over $F$ and such that the
resulting conjugation action of $\Gal(K|F)$ on $\Gal(D|K)$ is given by
$\rho$ will be understood in terms of the $\Gal(L_n|F)$-modules
$L_n^\times\!/L_n^{\times p}$ in characteristic~$0$ and
$L_n^+\!/\wp(L_n^+)$ in characteristic~$p$, where $\wp(x)=x^p-x$.
Therefore we study the filtered galoisian module $L^\times\!/L^{\times
p}$ (resp.~$L^+\!/\wp(L^+)$) for any finite tamely ramified split
galoisian extension $L$ of $F$ in the third
paper \citer\littlegalmods(), mostly following Iwasawa in
characteristic~$0$ and extending the results to characteristic~$p$.

Using these three ingredients, the set of primitive extensions $E$ of
the $p$-field $F$ is parametrised in the fourth and final
paper \citer\wildprim() of this series.  We show there that not only
the discriminant of $E$ but the filtered group $\Gal(\hat E|F)$ (where
$\hat E$ is the galoisian closure of $E$ over $F$) can be recovered
from its parameter, and illustrate the theory in the simplest case of
primitive quartic extensions of dyadic fields (or $2$-fields).  In
particular, the primitive quartic extensions of $\Q_2$ $\oldstyle(17)$
will be recovered without every having to write down a polynomial.

We also make some historical remarks in \citer\wildprim()~; suffice it
to say here that this method is a vast generalisation from the case
$n=1$ treated in \citer\monatshefte() and that similar but somewhat
less precise results have recently been obtained by Del Corso,
Dvornicich and Monge \citer\deldvomonge() when the $p$-field $F$ has
characteristic~$0$.


\smallbreak

\numeroter {\bf Sources}.  The proof of Galois's
theorem~$\oldstyle(3)$ is taken from Cox \citer\cox() and Neumann
\citer\neumann()~; both authors include detailed historical analyses.
The related statement $\oldstyle(4)$ is taken from \citer\cox().  Prompt and
elegant answers to the MathOverflow question \citer\csdmo()
were graciously provided by Geoff Robinson, Jeremy Rickard and Peter
Neumann.  I want to express here my heartfelt gratitude to them all.

\bigbreak
\unvbox\bibbox

\bye